\documentclass{article}
\usepackage{graphicx}
\usepackage{amsmath}
\usepackage{amsfonts}
\usepackage{amssymb}
\newtheorem{theorem}{Theorem}

\newtheorem{corollary}[theorem]{Corollary}

\newtheorem{definition}[theorem]{Definition}

\newtheorem{lemma}[theorem]{Lemma}

\newtheorem{proposition}[theorem]{Proposition}

\begin{document}
\subsubsection{\bigskip GAUSS INTEGERS AND DIOPHANTINE FIGURES}

Stancho Dimiev and Krassimir Markov

\ 

INTRODUCTION

The study of Diophantine figures in the plane (see [DT], [BDM], [D], [B])
involves different geometric and number-theoretic motions. First this is
Diophantine plane, i.e. the set of all points in Descartes plane
$\mathbb{R\times R}$ with integer coordinates . So Diophantine plane is just
the Cartesian product $\mathbb{Z\times Z}$. Z is the ring of all integers. In
terms of coordinate's operation in $\mathbb{Z\times Z}$ we develop a kind of
geometry of Diophantine figures, for instance constructions of some triangles
etc. However it is more reasonable to replace Descartes plane by Gauss plane
of all complex numbers $\mathbb{C=\{}a+ib:$ $a,b\in\mathbb{R}\},$ and
Diophantine plane by Gauss-Diophantine plane $\mathbb{Z+}i\mathbb{Z}%
=\{n+im:n,m\in\mathbb{R}\}.$ In fact $\mathbb{Z+}i\mathbb{Z}$ is defined as a
tensor product of $\mathbb{Z}$ with itself under the base $(1,0)$ and $(0,1).$
Gauss-Diophantine plane is stable under the operation ''multiplication'' which
give the possibility to define indecomposable elements. So Gauss-Diophantine
plane is a commutative ring without zero divisors, denoted ordinary by
$\mathbb{Z}[i].$ The elements of $\mathbb{Z}[i]$ are called Gauss integers.
This paper gives a survey on the obtained up to now results and proposes some
new results and problems.

\subsubsection{GAUSS INTEGERS}

\paragraph{GENERAL PROPERTIES, UNITIES, ASSOCIATE ELEMENT}

It is well known that the ring $\mathbb{Z}[i]$ is an Euclidean domain with
respect to the norm $\mathbf{N}(n+im)=$ $n^{2}+m^{2}.$ The norm is a function
of the type $\mathbb{Z}[i]\rightarrow\mathbb{N}$, such that for every two
elements $a$ and $b$ in $\mathbb{Z}[i],$ with $\ b$ different from zero, there
are two elements $c$ and $d$ in $\mathbb{Z}[i]$ for which $a=cb+d$ and in the
case $d\neq0$ we have $\mathbf{N(}d)<\mathbf{N(}b).$

An element $a+ib$ is an unit in $\mathbb{Z}[i]$ if there is $x+iy\in
\mathbb{Z}[i]$ such that
\[
(a+ib)(x+iy)=1,\text{ \ \ \ \ }1=1+i.0
\]

It is not difficult to see that there are only four units in $\mathbb{Z}[i],$
namely $1,$ $-1,i,-i.$

We say that $\beta\in\mathbb{Z}[i]$ is an associate element of $\alpha
\in\mathbb{Z}[i]$ if $\beta=\varepsilon\alpha,$ where $\varepsilon$ is an unit
in $\mathbb{Z}[i].$ Clearly if $\beta$ is an associate to $\alpha,$ then
$\alpha$ is an associate to $\beta.$ Every $\alpha\in\mathbb{Z}[i]$ has four
associate elements $1.\alpha,(-1).\alpha,i.\alpha,(-i).\alpha.$

\paragraph{PARITY IN $\mathbb{Z}[i]:$EVEN AND ODD GAUSS INTEGERS}

\begin{definition}
We say that the Gauss integers $n+im,$ $n,m\in\mathbb{Z},$ is an even integer
in $\mathbb{Z}[i]$ if \ $n$ and $m$ are of the same parity in $\mathbb{Z}.$ If
$n$ and $m$ are of different parity in $\mathbb{Z}$ we say that $n+im$ is an
odd Gauss integer in $\mathbb{Z}[i].$

\begin{proposition}
\textbf{1}. If the Gauss integer $n+im$ is even in $\mathbb{Z\ }[i]$ there is
$p+iq\in\mathbb{Z\ }[i]$ such that
\[
n+im=2(p+iq)\text{ \ \ or \ \ \ }n+im=2(p+iq)+1+i
\]
Respectively, if $n+im$ is odd in $\mathbb{Z}[i]$ there is $r+is$ such that
\[
n+im=2(r+is)+1\ \ \ \text{or}\ \ \ n+im=2(r+is)+i
\]
\end{proposition}

\begin{proposition}
\textbf{2}. The Gauss integer $n+im$ is even in $\mathbb{Z\ }[i]$ iff
\ $n+im\equiv0(\operatorname{mod}(1+i)),$ respectively $n+im$ is odd in
$\mathbb{Z\ }[i]$ iff \ $n+im\equiv1(\operatorname{mod}(1+i)).$
\end{proposition}

\begin{proposition}
\textbf{3. }The Gauss integer $n+im$ is even in $\mathbb{Z\ }[i]$ iff the norm
$\mathbf{N}(n+im)=$ $n^{2}+m^{2}$ is even in $\mathbb{Z}.$ Respectively,
$n+im$ is odd in $\mathbb{Z\ }[i]$ iff the norm $\mathbf{N}(n+im)$ is odd \ in
$\mathbb{Z}.$
\end{proposition}
\end{definition}

It is not difficult to prove that the above three propositions are equivalent.

\paragraph{ARITHMETIC PROPERTIES OF GAUSS INTEGERS}

It is easy to prove (with the help of the above stated propositions) that :

\textbf{The sum of Gauss integers }of common parity is even, and the sum of
Gauss integers of different parity is odd.

\textbf{The product of Gauss integers }satisfy the same rules as in
$\mathbb{Z}.$

\qquad$(even)\times(even)=(even)$, \ \ $(odd)\times(odd)=(odd),$
\ \ $(even)\times(odd)=(een).$

\textbf{The square }of an even Gauss integer is an even Gauss integer,
respectively the square of an odd Gauss integer is an odd Gauss integer. More
precisely, if $\alpha$ is even in $\mathbb{Z}[i]$ then $\alpha=(1+i)^{2}%
\alpha_{2}=4\alpha_{1}.$ If $\beta$ is odd in $\mathbb{Z}[i]$ then $\beta
^{2}=(1+i)^{2}\beta_{2}+1=4\beta_{1}+1.$

\textbf{Square radical }of a Gauss integer (S. Dimiev)

Let $\alpha=m+im$ be a Gauss integer with $m\neq0.$ We consider the equation
$z^{2}=\alpha$ with $z\in\mathbb{C}.$ Each solution of this equation is called
square radical of $\alpha.$ In the case that there is $l\in\mathbb{N}$ such
that $(n,m,l)$ to be a Pythagorean triple in $\mathbb{Z},$ i.e. $n^{2}%
+m^{2}=l^{2},$ we have
\[
\sqrt{\frac{n+im}{2}}=\frac{t+i\frac{m}{t}}{2},\text{ \ \ }t^{2}=n+l
\]

If we set $l-n=k^{2},$ then $m=kt,$ and the above formula will see as follows
\[
2\sqrt{\frac{n+im}{2}}=t+ik
\]

Proof (according to P.Guncheva).

We set $z=x+iy,$ $x,y\in\mathbb{R}.$ It follows that $x^{2}-y^{2}=n$ and
$2xy=m.$ Then we have $x\neq0$ and $y=\frac{m}{2x}$. The following biquadratic
equation for $x$ holds:%
\[
4x^{4}-4nx^{2}-m^{2}=0
\]
Setting $x^{2}=v$ we obtain:%
\[
v_{1}=%
\frac12
(n+\sqrt{n^{2}+m^{2}})\text{ \ and \ }v_{2}=%
\frac12
(n-\sqrt{n^{2}+m^{2}})
\]
As $v_{2}<0$ the equation $x^{2}=v_{2}$ has no real solutions. So%
\[
x^{2}=%
\frac12
(n+\sqrt{n^{2}+m^{2}})
\]
As by the condition $(n,m,l)$ is a Pyhtagorean triple in $\mathbb{Z},$ we
have:%
\[
x^{2}=%
\frac12
(n+l)\text{ \ \ and \ \ }y=\frac{m}{2x}%
\]
Now we set $n+l=t^{2}$ and obtain the above written formula.$\blacksquare$

\paragraph{INDECOMPOSABLE GAUSS INTEGERS}

\begin{definition}
Definition[IR] An element $\alpha$ of $\mathbb{Z\ }[i]$ is called
\textbf{indecomposable} Gauss integer or\textbf{ prime} Gauss integer if it is
impossible to present it as product of two elements $\lambda,\mu
\in\mathbb{Z\ }[i],$ both of which are not units, i.e. different from $1,-1,i,-i.$
\end{definition}

Below we shall give examples of Gauss prime integers.

The norm of a Gauss integer is sum of two squares, i.e. if $\alpha=n+im,$
$n,m\in\mathbb{Z},$ then $\mathbf{N}(\alpha)=n^{2}+m^{2}.$ With this in mind
we set $\mathbf{N}(\alpha)\equiv t(\operatorname{mod}4),$ where $t=0,1,2,3.$
It is not difficult to see that the case $t=3$ is impossible. More precisely,
we have:%
\[
\mathbf{N}(\alpha)=n^{2}+m^{2}\neq4s+3,\text{ \ \ }s\in\mathbb{Z}%
\]

So in the case $\mathbf{N}(\alpha)$ is odd, it follows $\mathbf{N}%
(\alpha)=4s+1,$ or $\mathbf{N}(\alpha)\equiv1(\operatorname{mod}4).$ In the
next exposition we need two well known theorems from the number theory in
$\mathbb{Z}.$

The first one is a theorem of Fermat (see Edwards [ \ ]): \textit{each prime
number in }$\mathbb{Z}$ \textit{such that }$\mathit{p}\equiv
1(\operatorname{mod}4)$ \textit{can be written as sum of two squares, i.e.
there exist }$a,b\in\mathbb{Z}$ \textit{such that }$p=a^{2}+b^{2}.$
\textit{When }$a,b$ \textit{are positive integers} $a$ \textit{is odd}, $b$
\textit{is even, the above mentioned representation is unique}.

The second theorem asserts: \textit{if} $(x,y)=1$ \textit{then} $x^{2}+y^{2}$
\textit{has at least one prime divisor of the form} $4k+1$ (see for instance
T. Tonkov, Figural numbers, Sofia, 1971, pp. 56-102).

Now we turn to the description of the prime Gauss integers. The mentioned
description is based on the comparison with prime integers in $\mathbb{Z}.$ We
will see that some prime integers in $\mathbb{Z}$ are decomposable in
$\mathbb{Z}[i]$ and that all prime Gauss integers are divisors of prime
integers in $\mathbb{Z}.$

More precisely we have the following (sugested by [IR])

\begin{lemma}
\textbf{1}. If $\alpha$ is a prime element of $\mathbb{Z\ }[i]$, then there is
a prime $p\in\mathbb{Z}$ such that $\alpha|$ $p.$
\end{lemma}

For the proof it is enough to write $\mathbf{N(}\alpha)=\alpha\overline
{a}=p_{1}p_{2}...p_{s}$ where $p_{k}$ are primes in $\mathbb{Z}.\blacksquare$

\begin{lemma}
\textbf{2}. If $\alpha$ is a Gauss integer and $\mathbf{N(}\alpha)$ is a prime
integer in $\mathbb{Z},$ then $\alpha$ is prime Gauss integer.
\end{lemma}

For the proof it is enough to remark that if $\alpha=\lambda\mu,$ with
$\lambda$ and $\mu$ non units, then $\mathbf{N(}\alpha)=\mathbf{N(}%
\lambda)\mathbf{N(}\mu).$ As $\mathbf{N(}\alpha)$ is a prime in $\mathbb{Z}$
it follows that $\mathbf{N(}\lambda)$ or $\mathbf{N(}\mu)$ is equal to 1,
which means that either $\lambda$ or $\mu$ is unit.$\blacksquare$

Based on the above proved two lemmas we give examples of prime Gauss integer
and decomposable in $\mathbb{Z}[i]$ prime integers in $\mathbb{Z}.$ These are
$1+i$ in $\mathbb{Z}[i]$ and 2 in $\mathbb{Z}.$ Indeed $\mathbf{N}(1+i)=2$ and
we apply Lemma 2; $2=(=i)(1+i)^{2}$ and applying Lemma 1 we obtain $(1+i)|2,$
$(-i)$ is an unit.$\blacksquare$

The exact description of the prime Gauss integers is given by the following

\begin{theorem}
The indecomposable elements or the prime integers in $\mathbb{Z\ }[i]$ are the following:
\end{theorem}

\textbf{a)} All prime integers in $\mathbb{Z}$ of the form $4k+3$ and all
their associate elements in $\mathbb{Z}[i];$

\textbf{b)} The number $1+i$ and its associates;

\textbf{c)} \ If $p$ is a prime integer in $\mathbb{Z}$ of the form $4k+1$ and
$\mathbf{N}(\alpha)=p,$ $\alpha=x+iy,$ $x,y\in\mathbb{Z},$ i.e. $p=x^{2}%
+y^{2},$ then $\alpha$ and $\overline{a}$ are indecomposable elements in
$\mathbb{Z}[i]$ with all their associates.

\textbf{d)} There are no other indecomposable elements in $\mathbb{Z}[i]$,
more precisely if $\mathbf{N}(\alpha)=q$, where $q$ is not prime integer in
$\mathbb{Z}$, then $\alpha$ is indecomposable element in $\mathbb{Z}[i].$

Proof. \textbf{a). }Let $s\in\mathbb{Z}$ be a prime number in $\mathbb{Z}$ of
the form $s=4k+3.$ As $\mathbf{N}(s)=s^{2},$ supposing the validity of the
decomposition $s=\lambda\mu,$ $\lambda=a+ib,$ $\mu=c+id,$ we obtain
$s^{2}=(a^{2}+b^{2})(c^{2}+d^{2}),$ which is impossible. Indeed $s\neq
(a^{2}+b^{2}),(c^{2}+d^{2})$ in view of $s=4k+3$ (see a remark above). Of
course, $s^{2}=(a^{2}+b^{2})(c^{2}+d^{2}).$ This asserts that $s$ is
indecomposable in $\mathbb{Z}[i].$

\textbf{b). }In this case the proof was given above.

\textbf{c). }Let $\alpha=x+iy$ be non-trivial decomposable in $\mathbb{Z}[i]$,
say $\alpha=\lambda\mu,\lambda=a+ib,$ $\mu=c+id,\lambda$ and $\mu$ are both
non-units in $\mathbb{Z}[i].$ It follows that $p=x^{2}+y^{2}=(a^{2}%
+b^{2})(c^{2}+d^{2}).$ In view that $p$ is prime number in $\mathbb{Z}$, we
obtain that either of $(a^{2}+b^{2})$ or $(c^{2}+d^{2})$ must be 1, which
implies that it is an unit. This means that the admitted decomposition is
trivial, which is a contradiction.

\textbf{d).} If $\alpha=x+iy,q=x^{2}+y^{2}.$ We can suppose that $(x,y)=1.$ As
$q$ is sum of two squares, according to the above cited theorem $q=pt,$ with
$t>1$ and $p$ is prime in $\mathbb{Z}$ of the form $4k+1.$ Applying the cited
theorem of Fermat we obtain $p=u^{2}+v^{2}.$ Now setting $\beta=u+iv$ we
receive $x^{2}+y^{2}=(u^{2}+v^{2})t.$ According to \textbf{c). }we have that
$\beta$ and $\overline{\beta}$ are prime Gauss integers. (indecomposable
elements in $\mathbb{Z}[i])$. Clearly $\alpha\overline{\alpha}=\beta
\overline{\beta}t,$ and , for instance, $\alpha=\beta t_{1},$ $\overline
{\alpha}=\overline{\beta}t_{2},$ where $t_{1}t_{2}=t.$ As $t>1$ at least
$t_{1}$ or $t_{2}$ must be greater then 1. It gives a nontrivial decomposition
of $\alpha.\blacksquare$

\paragraph{GAUSS-PYTHAGOREAN INTEGERS}

\begin{definition}
(S. Dimiev) The Gauss integer $\alpha=x+iy,$ $x,y\in\mathbb{Z},$ is said to be
a Gauss-Pythagorean number if there exists $z\in\mathbb{Z}$ such that the
triple $x,y,z$ be a Pythagorean triple: $x^{2}+y^{2}=z^{2}.$
\end{definition}

We shall denote the set of all Gauss-Pythagorean integers by $\mathbb{GP}[i],$
$\mathbb{GP}[i]\subset\mathbb{Z}[i].$ The zero element and the units of
$\mathbb{Z}[i]$ are not Gauss-Pythagorean integers. We remark that the sum of
two Gauss-Pythagorean numbers is not Gauss-Pythagorean number in general, but
it is easy to prove that the product of two Gauss-Pythagorean numbers is
always a Gauss- Pythagorean integer. So $\mathbb{GP}[i]$ is a multiplicative
subsemigroup of the multiplicative group $\mathbb{Z}[i]-\{0\}$ of the ring
$\mathbb{Z}[i].$

The conjugate and the associate elements of a Gauss-Pythagorean element are
Gauss-Pythagorean integers too. The following Lemma is useful in the following exposition

\begin{lemma}
(K. Markov) Let $\alpha=x+iy$ be a Gauss-Pythagorean number them there exists
a Gauss integer $\tau$ such that $\mathbf{N(\alpha)=(N}(\tau))^{2}.$
\end{lemma}

Proof: First let we assume that $(x,y)=1.$ As $\alpha\in\mathbb{GP}[i]$ there
exist $z\in\mathbb{N}$ for which $(x,y,z)$ is a primitive Pythagorean triple
in $\mathbb{Z.}$ According to the well known representation of Pythagorean
triples in $\mathbb{Z},$ there are $p,q\in\mathbb{N}$such that $z=p^{2}%
+q^{2}.$ We set $\tau=p+iq$ which gives the statement. In the general case
$(x,y)=d>1,$ we have $\mathbf{N}(\alpha)=d^{2}\mathbf{N}(x_{1}+iy_{1})$ with
$x_{1},y_{1}\in\mathbb{Z}$ and $(x_{1},y_{1})=1.$ Setting $\tau=$
$x_{1}+iy_{1}$ we apply the first remark above: $\mathbf{N}(\alpha
)=d^{2}(\mathbf{N}(\tau_{1}))^{2}=\mathbf{N}(d\tau_{1}),$ with $d\tau_{1}%
\in\mathbb{Z}[i].\blacksquare$

\begin{corollary}
There is no element in $\mathbb{GP}$\/$\ [i]$ which is prime Gauss integer.
\end{corollary}

\begin{definition}
An element $\tau\in\mathbb{GP\ }[i]$ is said to be prime Gauss-Pythagorean
integer if is impossible to represent it as product of two elements of
$\mathbb{GP\ }[i].$
\end{definition}

\begin{theorem}
There exists an infinite number of indecomposable Gauss-Pythagorean numbers.
\end{theorem}

Proof: (K. Markov) Let we remark that if $\alpha\in\mathbb{GP}[i]$ and
$\mathbf{N}(\alpha)=p^{2},$ where $p$ is prime integer in $\mathbb{Z}$, then
$\alpha$ is prime integer in $\mathbb{GP}[i].$ Indeed supposing $\alpha
=\beta\gamma$ with $\beta,\gamma\in\mathbb{GP}[i],$ $\mathbf{N}(\beta)>1,$
$\mathbf{N}(\gamma)>1,$ we apply the Lemma proved above and obtain:
$\mathbf{N}(\beta)=(\mathbf{N}(\beta_{1}))^{2}$ and $\mathbf{N}(\gamma
)=(\mathbf{N}(\gamma_{1}))^{2}.$ Thus in $\mathbb{Z}$ we have $p^{2}%
=\mathbf{N}(\beta)\mathbf{N}(\gamma)=(\mathbf{N}(\beta_{1})\mathbf{N}%
(\gamma_{1}))^{2},$ $\beta_{1},\gamma_{1}\in\mathbb{GP}[i],$ which implies
$p=\mathbf{N}(\beta_{1})\mathbf{N}(\gamma_{1})$ with $\mathbf{N}(\beta
_{1})>1,$ $\mathbf{N}(\gamma_{1})>1.$ But this contradicts the condition $p$
to be prime integer in $\mathbb{Z}.$ So $\alpha$ is prime integer in
$\mathbb{GP}[i].$

Now we take a prime integer $p=4k+1.$ By the above cited Fermat theorem there
are natural integers $t$ and $s$ such that $p=t^{2}+s^{2}.$ We shall consider
the Gauss integer $t+is.$ Having in mind that $(t+is)^{2}=t^{2}-s^{2}+2its,$
we consider the triple $(2ts,t^{2}-s^{2},t^{2}+s^{2}).$ This triple is a
Pythagorean triple according to the well known formulae for natural
Pythagorean numbers. So $(t+is)^{2}$ is a Gauss-Pythagorean triple. The
mapping $p\mapsto(t+is)^{2}$ is an injective mapping between the sets: $\{p$
is prime : $p=4k+1\}\rightarrow\mathbb{GP}[i].$ Indeed if $p_{1}\neq p_{2}$
are two numbers from the first set then the corresponding $(t_{j}+is_{j}%
)^{2},$ $j=1,2$ are different because $\mathbf{N}((t_{j}+is_{j})^{2})=p_{j}$
are different for $j=1,2.$ Finally it is to be recalled that there are an
infinite number of prime integers in the progression $4k+1.\blacksquare$

\paragraph{PRIMITIVE TRIPLES IN $\mathbb{Z}[i]$}

\begin{definition}
We say that the triple $\alpha,\beta,\gamma$ of Gauss integers is
\textbf{primitive triple} if the unique common divisors of the elements of the
triple are the unities in $\mathbb{Z\ }[i],$ i.e. $1,-1,i,-i.$
\end{definition}

We denote this by $(\alpha,\beta,\gamma)=1.$ Like in $\mathbb{Z}$ we shall
write $(\alpha,\beta)=1$ when $\alpha$ and $\beta$ satisfy the same condition.
Ordinary by $(\alpha,\beta,\gamma)$ is denoted the GCD of $\alpha,\beta$ and
$\gamma.$ Analogously $(\alpha,\beta)$ is the GCD of $\alpha$ and $\beta.$ It
is easy to see that $(\alpha,\beta,\gamma)=((\alpha,\beta),\gamma
)=(\alpha,(\beta,\gamma))=((\alpha,\gamma),\beta).$ If $\overline{\alpha}$ is
the complex conjugate of $\alpha$ we have: $(\alpha,\beta)=\delta$ iff
$(\overline{\alpha},\overline{\beta})=\overline{\delta}$ in particular
$(\alpha,\beta)=1$ iff $(\overline{\alpha},\overline{\beta})=1.$ Analogously
$(\alpha,\beta,\gamma)=\delta$ iff $(\overline{\alpha},\overline{\beta
},\overline{\gamma})=\overline{\delta}$ and in particular $(\alpha
,\beta,\gamma)=1$ iff $(\overline{\alpha},\overline{\beta},\overline{\gamma
})=1.$ Clearly $(\alpha,\beta)=(\alpha,\gamma)=(\beta,\gamma)=1$ implies
$(\alpha,\beta,\gamma)=1$ but the inverse is not true.

\begin{proposition}
\textbf{4}. $(\mathbf{N(}\alpha),\mathbf{N}(\beta),\mathbf{N(}\gamma))=1$ in
$\mathbb{Z}$ implies $(\alpha,\beta,\gamma)=1$ in $\mathbb{Z}[i].$
\end{proposition}

Proof: Let we remark that if $\lambda$ is a divisor of $\alpha,\beta$ or
$\gamma$ then $\mathbf{N}(\lambda)$ is a divisor of $\ \mathbf{N(}%
\alpha),\mathbf{N}(\beta)$ or $\mathbf{N(}\gamma).$ This follows from
$\alpha=\lambda\eta$ implies $\mathbf{N}(\alpha)=\mathbf{N}(\lambda
)\mathbf{N}(\eta).$

Now supposing that $(\alpha,\beta,\gamma)=\delta$ with $\mathbf{N}(\delta)>1$
we receive that $\mathbf{N}(\delta)$ is a divisor of $(\mathbf{N(}%
\alpha),\mathbf{N}(\beta),\mathbf{N(}\gamma))$ which is impossible. The
inverse is not true.

\textit{Example: }$(3+i,2+i,8+i)=1,$ but \textit{ }$(\mathbf{N}%
(3+i),\mathbf{N}(2+i),\mathbf{N}(8+i))=(10,5,65)=5.$ We remark that
$(3+i,2+i)=1$ does not imply $(\mathbf{N}(3+i),\mathbf{N}(2+i))=1.\blacksquare$

\paragraph{PYTHAGOREAN TRIPLES, PRIMITIVE PYTHAGOREAN TRIPLES}

In a primitive Pythagorean triple in $\mathbb{Z}[i]$ at least one element is
odd. Let $\alpha,\beta,\gamma$ be a primitive Pythagorean triple, i.e.
$\alpha^{2}+\beta^{2}=\gamma^{2}$ and $(\alpha,\beta,\gamma)=1.$ There are two
possibilities for $\gamma:$ the first one is $\gamma$ to be even. In this case
we have:%
\[
(odd)^{2}+(odd)^{2}=(even)^{2}%
\]

The second possibility is $\gamma$ to be odd. So we have for instance:%
\[
(even)^{2}+(odd)^{2}=(odd)^{2}%
\]

Proposition 4\textbf{ }obtains for Pythagorean triples the following stronger form:

\begin{proposition}
\begin{proposition}
\textbf{4'. } $(\mathbf{N(}\alpha),\mathbf{N}(\beta),\mathbf{N(}\gamma))=1$ in
$\mathbb{Z}$ is equivalent to $(\alpha,\beta,\gamma)=1$ in $\mathbb{Z\ }[i].$
\end{proposition}
\end{proposition}

Proof: It is sufficient to prove that $(\alpha,\beta,\gamma)=1$ in
$\mathbb{Z}[i]$ implies $(\mathbf{N(}\alpha),\mathbf{N}(\beta),\mathbf{N(}%
\gamma))=1$ in $\mathbb{Z}.$ Let $(\mathbf{N(}\alpha),\mathbf{N}%
(\beta),\mathbf{N(}\gamma))=d$ in $\mathbb{Z}.$Let $p$ e a prime divisor of
$d.$ Clearly $d|\mathbf{N}(\alpha),d|\mathbf{N}(\beta),d|\mathbf{N}(\gamma).$

Let $\mathbf{N}(\alpha)=\alpha\overline{\alpha}=\alpha_{1}d,\mathbf{N}%
(\beta)=\beta\overline{\beta}=\beta_{1}d,\mathbf{N}(\gamma)=\gamma
\overline{\gamma}=\gamma_{1}d.$ $p|d$ implies: $p|\alpha$ or $p|\overline
{\alpha},$ $p|\beta$ or $p|\overline{\beta},$ $p|\gamma$ or $p|\overline
{\gamma}.$ If \ $p|\alpha$ and $p|\beta$ it follows \ $p|\gamma$ and this
contradicts the condition $(\alpha,\beta,\gamma)=1.$ If $p|\overline{\alpha}$
and $p|\overline{\beta}$ it follows \ $p|\overline{\gamma}$ and contradicts
the condition $(\overline{\alpha},\overline{\beta},\overline{\gamma
})=1.\blacksquare$

\begin{proposition}
\textbf{5. }Formulae for primitive Pythagorean triples:

$\alpha=2\lambda\mu,$ $\beta=\lambda^{2}-\mu^{2},\gamma=\lambda^{2}+\mu
^{2},(\lambda,\mu)=1,\lambda,\mu\in\mathbb{Z}\ [i]$

where $\lambda$ and $\mu$ are of different parity.
\end{proposition}

Proof: By condition $\alpha^{2}+\beta^{2}=\gamma^{2}$ and $(\alpha
,\beta,\gamma)=1.$ As $\alpha^{2}=\gamma^{2}-\beta^{2}=(\gamma-\beta
)(\gamma+\beta),$ let $\delta=(\gamma-\beta,\gamma+\beta).$ It follows that
$\gamma-\beta=\delta\gamma_{1},\gamma+\beta=\delta\gamma_{2}$ and also
$2\gamma=\delta(\gamma_{2}+\gamma_{1}),2\beta=\delta(\gamma_{2}-\gamma_{1}).$
So $\delta$ is a divisor of $(2\gamma,2\beta)=$ $2(\gamma,\beta)=2$ as
$(\gamma,\beta)=1.$ It is proved that $\delta$ is a divisor of $2.$ The
divisors of $2$ in $\mathbb{Z}[i]$ are $1,1+i,2$ up to an associate element in
$\mathbb{Z}[i].$

In the case $\delta=1,$ having in mind that $\gamma$ and $\beta$ are of common
parity, we conclude that $\gamma-\beta$ and $\gamma+\beta$ are both even and
we can write $\gamma-\beta=(1+i)\delta_{1},\gamma+\beta=(1+i)\delta_{2}.$ But
this means that $1+i$ is a divisor of $1=(\gamma-\beta)(\gamma+\beta),$ which
is impossible.

In the case $\delta=1+i$ we have $\gamma-\beta=(1+i)\gamma_{1}$ and
$\gamma+\beta=(1+i)\gamma_{2}$ and, as corollary, $1+i=(1+i)(\gamma_{1,}%
\gamma_{2})$ or $(\gamma_{1,}\gamma_{2})=1.$ So we obtain $\alpha
^{2}=(1+i)^{2}\gamma_{1}\gamma_{2}$ and also $\left(  \frac{\alpha}%
{1+i}\right)  ^{2}=\gamma_{1}\gamma_{2}.$ According to a well known lemma
there are $\rho,\sigma$ such that $\gamma_{1}=\rho^{2}$ and \ $\gamma
_{2}=\sigma^{2}.$ This implies: $\alpha=(1+i)\rho\sigma,2\gamma=(1+i)(\rho
^{2}+\sigma^{2}),2\beta=(1+i)(\sigma^{2}-\rho^{2}).$ If $\rho$ and $\sigma$
are of common parity $\rho^{2}+\sigma^{2}$ and $\sigma^{2}-\rho^{2}$ are both
even . This implies that $1+i$ is a common multiple for $\alpha,\beta$ and
$\gamma$ in contradiction with the primitivity of the triple $\alpha
,\beta,\gamma.$ If $\rho$ and $\sigma$ are of different parity $\gamma$ and
$\beta$ will be not Gauss integers as $1/(1+i)$ is not a Gauss integer.

Finally we have the case $\delta=2.$ In this case we have $2=(\gamma
-\beta)(\gamma+\beta),$ and therefore $(\gamma_{1,}\gamma_{2})=1.$ Applying
the same Lemma for $\left(  \frac{\alpha}{2}\right)  ^{2}=\gamma_{1}\gamma
_{2}$ we conclude that the are $\tau$ and $\kappa$ in $\mathbb{Z}[i]$ such
that $\gamma_{1}=\tau^{2},\gamma_{2}=\kappa^{2},(\tau,\kappa)=1.$ It follows:
$\alpha=2\tau\kappa,\beta=\gamma_{2}-\gamma_{1}=\kappa^{2}-\tau^{2}%
,\gamma=\kappa^{2}+\tau^{2}.$ In view that $(\alpha,\beta,\gamma)=1,$ then
$\tau$ and $\kappa$ must be of different parity.$\blacksquare$

\begin{proposition}
\textbf{6.}(Fermat Last Theorem in $\mathbb{Z\ }[i]$ for $n=4).$ The equation
$x^{4}+y^{4}=z^{4}$ has no integer solutions.
\end{proposition}

Proof: First let consider the equation $x^{4}+y^{4}=z^{2}.$ Supposing that
there exists in $\mathbb{Z}[i]$ a solution $\alpha,\beta,\gamma$, we apply the
Proposition 5 to the triple $(\alpha^{2},\beta^{2},\gamma)$ which will be a
solution of the equation $x^{2}+y^{2}=z^{2}.$ According to Proposition 5 there
are two Gauss integers $\kappa$ and $\tau$, such that:%
\[
\alpha^{2}=2\kappa\tau,\beta^{2}=\kappa^{2}-\tau^{2},\gamma=\kappa^{2}%
+\tau^{2}%
\]

It is easy to prove that the combination $\kappa=even$ and $\tau=odd$ is
impossible. Thus we have $\kappa=odd$ and $\tau=even.$ Let $\tau=2\xi.$ We
receive $%
\genfrac{(}{)}{}{}{\alpha}{2}%
^{2}=\kappa\tau$ and we apply the Lemma to the equation $W^{2}=UV$ in
$\mathbb{Z}[i].$ There are two Gauss integers $\kappa_{1}$ and $\xi_{1},$ such
that: $\kappa=\kappa_{1}^{2},\xi=\xi_{1}^{2},(\kappa_{1},\xi_{1})=1.$ It
follows that $\tau=2\xi_{1}$ and $(2\xi_{1}^{2})^{2}+\beta^{2}=(\kappa_{1}%
^{2})^{2}.$ Applying again Proposition 5:%
\[
2\xi_{1}^{2}=2\eta\zeta,\beta^{2}=\eta^{2}-\zeta^{2},\kappa_{1}^{2}=\eta
^{2}+\zeta^{2}%
\]

Now we apply again the Lemma to the equation $W^{2}=UV$ in $\mathbb{Z}[i]$ to
the first equation above. There are two Gauss integers $\eta_{1},\zeta_{1}$
such that $\eta=\eta_{1}^{2},\zeta=\zeta_{1}^{2},(\eta_{1},\zeta_{1})=1.$
Replacing these two integers in the third equation we obtain:%
\[
\eta_{1}^{4}+\zeta_{1}^{4}=\kappa_{1}^{2}%
\]

We can choose $\gamma$ in the solution $\alpha,\beta,\gamma$ to be with
minimal norm $\mathbf{N}(\gamma).$ With this in mind we have: $\mathbf{N}%
(\gamma)<\mathbf{N}(\kappa_{1}).$ So%
\[
\mathbf{N}(\gamma)<\mathbf{N}(\kappa_{1})<(\mathbf{N}(\kappa_{1}%
))^{2}=\mathbf{N}(\kappa)
\]

On the other hand : $\mathbf{N}(\gamma)=\mathbf{N}(\kappa^{2}+\tau
^{2})=(\mathbf{N}(\kappa))^{2}+(\mathbf{N}(\tau))^{2}+2\operatorname{Re}%
(\kappa\overline{\tau})^{2}.$ In the case $\operatorname{Re}(\kappa
\overline{\tau})^{2}\geq0$ we obtain:%
\[
(\mathbf{N}(\kappa))^{2}+(\mathbf{N}(\tau))^{2}+2\operatorname{Re}%
(\kappa\overline{\tau})^{2}\leq\mathbf{N}(\kappa)
\]

which is impossible as $\mathbf{N}(\kappa)$ is a natural integer $\geq1.$

\textit{In the case }$\operatorname{Re}(\kappa\overline{\tau})^{2}\leq0$ \ the
proposed method of a Fermat type desante with respect to the norm  falls. A
more sophisticated  method is developted, but it is more longer and will be
exposed elsewhere. .

\paragraph{DIOPHANTINE FIGURES}

We shall consider the Cartesian plane $\mathbb{R\times R}$ ($\mathbb{R}-$ the
field of real numbers). A \textbf{complete Cartesian graph }is by
definition\textbf{ }the couple (\textbf{V,S}) where \textbf{V }is the set of
points in $\mathbb{R\times R},$ called vertices, and \textbf{S} is the set of
all segments $[P.Q]$ with $P,Q\in\mathbf{V}$, $P\neq Q.$ A \textbf{Cartesian
Erd\"{o}s graph} is by definition a Cartesian graph $(\mathbf{V,S})$ for which
the length of each segment in $\mathbf{S}$ is a integer number $\neq0.$ If the
set of vertices $\mathbf{V}$ is infinite we shall say that $(\mathbf{V,S})$ is
an infinite graph.

\begin{theorem}
(Erd\"{o}s): The vertices of an infinite Cartesian Erd\"{o}s graph are
situated on a straight line in the Cartesian plane.

\begin{definition}
The Cartesian product $\mathbb{Z\times Z}$ will be called \textbf{Diophantine
plane}.
\end{definition}
\end{theorem}

Clearly $\mathbb{Z\times Z\subset R\times R},$ or the Diophantine plane is the
lattice of the points in $\mathbb{R\times R}$ with integer coordinates. We
will consider complete graphs in the Diophantine plane, i.e. the set of
couples $(\mathbf{V,S}),$ $\mathbf{V\subset\ }\mathbb{Z\times Z}$ and
$\mathbf{S}$ is the same as above. A \textbf{Diophantine figure }is by
definition a complete graph in the Diophantine plane for which the length of
each of its segments is an integer number. Diophantine figures which contain
at least three different non-collinear vertices will be considered.
\textbf{Erd\"{o}s-Diophantine figure }is by definition a\textbf{ maximal
}Diophantine figure, i.e. a Diophantine figure for which there is no a larger one.

The existence of Erd\"{o}s-Diophantine figures follows from the above cited
Erd\"{o}s theorem. Indeed, according this theorem each increasing sequence of
finite non-linear Diophantine figures $F_{1}\subset F_{2}\subset...\subset
F_{n}\subset$ ...stabilizes at some index $k_{0}\in\mathbb{N}.$ Then
$F_{k_{0}}$ is an Erd\"{o}s-Diophantine figure.

\textbf{Diophantine planimetry: examples}

A closed path in a Diophantine figure $F$ is defined by a sequence of vertices
$P_{0},P_{1},\ldots,P_{n}=P_{0},$ $P_{j}\in F,$ and the union of the
connecting segments. For a Diophantine triangle there is only one closed path
constituted from all vertices and all segments of the triangle.

\begin{proposition}
\textbf{6. }The sum of \ lengths of the segments of a closed path in a
Diophantine figure is an even integer.
\end{proposition}

The proof can be derived by induction from the following

\begin{lemma}
(M. Brancheva) The sum of lengths of the sides of a Diophantine triangle is
always an even integer.$\blacksquare$
\end{lemma}

This Lemma is generalization of the analogous property of Pythagorean triangles.

\begin{proposition}
\textbf{7. }(M. Brancheva) Let $\triangle ABX$ \ lies in the Diophantine plane
and $(a_{1},a_{2}),(b_{1},b_{2}),(x_{1},x_{2})$ are the coordinates of the
vertices $A,B,X$ resp. Let suppose $a_{1}=a_{2}=0;$ let the lengths of the
segments of the triangle are: $|AB|=\mathbf{c,|}BC\mathbf{|=a}%
,|AC|=\mathbf{b.}$ For given lengths $\mathbf{a,b,c}$ and given coordinates
$(b_{1},b_{2})$ we have the following Diophantine equation of first degree for
$x_{1},x_{2}:$%
\[
2b_{1}x_{1}+2b_{2}x_{2}=\mathbf{b}^{2}+\mathbf{c}^{2}-\mathbf{a}^{2}%
\]
\end{proposition}

Proof: Clearly we have $\mathbf{c}^{2}=b_{1}^{2}+b_{2}^{2}$ , $\mathbf{b}%
^{2}=x_{1}^{2}+x_{2}^{2}.$ On the other hand we have $\mathit{a}^{2}%
=(b_{1}-x_{1})^{2}+(b_{2}-x_{2})^{2}.$ The last three equalities imply the
statement of the proposition. $\blacksquare$

According to the theory of Diophantine equations of first degree if
$(x_{1}^{\ast},x_{2}^{\ast})$ is a solution of the equation then all solutions
are given by the formulae:%
\[
x_{1}=x_{1}^{\ast}+\frac{b_{2}t}{(b_{1},b_{2})},\text{ \ }x_{2}=x_{2}^{\ast
}-\frac{b_{1}t}{(b_{1},b_{2})}%
\]

where $t\in\mathbb{Z}$ and $(b_{1},b_{2})$ is the GCD of $b_{1},b_{2}.$
Eliminating $t$ from the above two equations we obtain for the solutions
$x_{1},x_{2}:$%
\[
x_{2}-x_{2}^{\ast}=-\frac{b_{1}}{b_{2}}(x_{1}-x_{1}^{\ast})
\]

Thus the point with coordinates $(x_{1},x_{2})$ lies on the perpendicular to
the segment $AB$ through the point with coordinates $(x_{1}^{\ast},x_{2}%
^{\ast}).$

Fig. 1

\textbf{Remark:}The above written Diophantine equation is not always solvable.
Indeed, if a solution exists the number $\mathbf{b}^{2}+\mathbf{c}%
^{2}-\mathbf{a}^{2}$ must be even, but according to Proposition 7 the same is
true for the number $\mathbf{b}^{2}+\mathbf{c}^{2}+\mathbf{a}^{2}.$ This
implies the following necessary condition: the number $\mathbf{b}%
^{2}+\mathbf{c}^{2}$ must be even. We see that the classical construction is
not always possible for Diophantine triangles.

\paragraph{DIOPHANTINE TRIANGLES: CLASSIFICATION}

Each Diophantine triangle can be inscribed in a uniquely determined rectangle
with sides parallel to the coordinate axes. with the help of this enveloping
rectangle we can formulate the following

\begin{lemma}
(Classification lemma): There are 4 essentially different types of
\ Diophantine triangles (see the Figures below):

(1) Pythagorean triangle;

(2) and (3) obtained from 2 Pythagorean triangles with common cathetus;

(4) new kind of Diophantine triangle.
\end{lemma}

The proof can be obtained by simple examination of the possibilities for the
disposition of the vertices of the inscribed triangle on the sides of the
enveloping rectangle.

Fig. 2

\textbf{Remark:}The supplementary part of the enveloping rectangle with
respect to the inscribed Diophantine triangle is composed by 1, 2 or 3
Pythagorean triangles.

The Classification lemma suggests the possibility of different calculations.
We shall consider only the case (4).

Applying the well known formulae for Pythagorean triples we can write:

for $\triangle APB:$ \ \ $AP=u^{2}-v^{2},BP=2uv,AB=u^{2}-v^{2},$ $\ \ u>v$

for $\triangle BQC:$ \ \ $CQ=p^{2}-q^{2},BQ=2pq,BC=p^{2}+q^{2},$ $\ p>q$

for $\triangle ACR:$ \ \ $RC=x^{2}-y^{2},AR=2xy,AC=x^{2}+y^{2},$ $\ \ x>y$

From $AR=PQ$ it follows: $2xy=2pq+2uv$ and from $AP=RQ$ it follows:
$u^{2}-v^{2}=x^{2}-y^{2}+p^{2}-q^{2}.$

In the Gauss-Diophantine plane we have:%
\[
(x+iy)^{2}=(q+ip)^{2}+(u+iv)^{2}%
\]

This means that $(q+ip,u+iv,x+iy)$ is a Pythagorean triple in $\mathbb{Z}[i].$
Now we apply the Proposition 5, according to which:

$x+iy=(a+ib)^{2}+(c+id)^{2},$ \ \ $a,b,c,d\in\mathbb{Z}$

$u+iv=(a+ib)^{2}-(c+id)^{2},$

$q+ip=2(a+ib)(c+id).$

Consequently:

$x=a^{2}-b^{2}+c^{2}-d^{2},y=2(ab-cd),u=a^{2}-b^{2}-c^{2}+d^{2},v=2(ab-cd),p=2(ad+bc),q=2(ac-bd).$

Taking $b=c=d=1$ we get only one parameter $a.$ After calculations we obtain
for $a=4$ the following triple: $AB=261,BC=136,AC=325.$ For $a=5$ we obtain
another triple:$AB=640,BC=208,AC=270.$ It can be verified that these two
triples define Diophantine triangles, and the supplementary triangles
$APB,BQC,CRA$ are Pythagorean ( $AP=189,PB=180,BQ=120,QC=64,CR=125,RA=300$ ).

\textbf{Remark: }The above exposed examples give an idea how to proceed
practically for to get Diophantine triangles of kind 4.

\paragraph{\bigskip DIOPHANTINE FIGURES COMPOSED BY PYTHAGOREAN TRIANGLES WITH
COMMON CATHETUS}

The simplest Diophantine figures of many vertices are composed by many
Pythagorean triangles (see Fig.3)

Fig. 3.

We shall consider the set of Pythagorean triples $(x,n,z)$ in $\mathbb{N},$ i.
e. $x^{2}+n^{2}=z^{2}.$ It is not supposed that these triples are primitive.
We introduce the function $\kappa:\mathbb{N\rightarrow N},$ $n\longmapsto
\kappa(n),$ where $\kappa(n)$ is the number of all Pythagorean triangles with
cathetus $n.$ By definition $\kappa(0)=0$ and $\kappa(m)=0$ if $m$ is not a
cathetus in a Pythagorean triangle.

We shall write $d$\/$\ |\ n$ if $d$ is a divisor of $n.$ By $\eta(d)$ is
denoted the set of all primitive Pythagorean triples with $d$ as cathetus. It
is clear that from each primitive Pythagorean triple $(u,d,v)$ we can obtain a
Pythagorean triple $(x,n,y)$ with cathetus $n.$ It is sufficient to multiply
by $\frac{n}{d}=k,$ i.e. $x=\frac{n}{d}u,$ $y=\frac{n}{d}v.$ Having in mind
all divisors of $n$ we can introduce the following formula:%
\[
\kappa(n)=\sum\limits_{d\/\ |\ n}\eta(d)
\]

\begin{lemma}
If $\delta(d)$ is the set of all divisors of $d$ then $\eta(d)<\delta(d).$
\end{lemma}

Proof: As we have a primitive Pythagorean triple $(u,d,v)$ we can use the
representation:%
\[
u=2st,d=s^{2}-t^{2},v=s^{2}+t^{2},\text{ \ \ }s,t\in\mathbb{N}\text{\ }%
\]

Then $d=(s-t)(s+t).$ If $p$ is a divisor of $d$ we can set $s-t=p$ and
$s+t=\frac{d}{p}.$ The number of solutions in $\mathbb{Z}$ of the above
written system may be $0$ in general. This shows that $\delta(d)>\eta(d).$

Now let $n$ be as follows: $n=p_{1}^{\alpha_{1}}p_{2}^{\alpha_{2}}%
...p_{r}^{\alpha_{r}}$ and \ $d=p_{1}^{\beta_{1}}p_{2}^{\beta_{2}}%
...p_{r}^{\beta_{r}}$ with $0\leq\beta_{j}\leq\alpha_{j}.$ We have:%
\begin{align*}
\delta(p_{1}^{\beta_{1}}p_{2}^{\beta_{2}}...p_{r}^{\beta_{r}})  &
=(1+\beta_{1})(1+\beta_{2})\cdots(1+\beta_{r})\\
\sum\limits_{d\/\ |\ n}\beta(d)  &  =\sum\limits_{d\/\ |\ n}(1+\beta
_{1})(1+\beta_{2})\cdots(1+\beta_{r})
\end{align*}

So we have:%
\[
\kappa(n)<\delta(n)<\left(  \prod\limits_{j=1}^{n}(1+\alpha_{j})\right)  ^{2}%
\]

\begin{theorem}
(M. Brancheva): \ $\kappa(n)=O(n^{\varepsilon}).$
\end{theorem}

Proof: According to a well known formula (see Prachar [P])%
\[
\delta(n)<\exp\{(1+\rho)\frac{\ln2\ln n}{\ln\ln(n)}\}
\]

But $\exp\{(1+\rho)\frac{\ln2\ln n}{\ln\ln(n)}\}=n^{2(1+\rho)\ln2/\ln\ln(n)}.$
We receive:%
\begin{align*}
\lbrack\delta(n)]^{2}  &  <n^{2(1+\rho)\ln2/\ln\ln(n)}\\
\lbrack\delta(n)]^{2}  &  <O(n^{\varepsilon}),\varepsilon>0.\blacksquare
\end{align*}

\textbf{Application:}$\underset{n\rightarrow\infty}{\lim}\frac{\kappa(n)}%
{n}=0$ which is conjecture of S. Dimiev.

It is proved in stronger form: \ \ $\underset{n\rightarrow\infty}{\lim}%
\frac{\kappa(n)}{n^{\varepsilon}}=0$ \ \ for every $\varepsilon>0.$

\subsection{\bigskip PROBLEMS}

\textbf{1.} Let us denote by $\chi(l)$ the number of all Pythagorean triangles
with hypotenuse $l,l\in\mathbb{N}.$ Find the asymptotic of the function
$\chi(l)$ when $l\rightarrow\infty$ following the above exposed case of the
function $\kappa(n).$

\textbf{2. }Given a Diophantine triangle $ABC$ is it possible to find a point
$D$ in the Diophantine plane such that $ABCD$ to be a Diophantine figure? The
case when there is no such point $D$ means that the triangle $ABC$ is an
Erd\"{o}s-Diophantine figure. Are there Erd\"{o}s-Diophantine triangles? In
the case when there is such a point $D$ is it possible to find an effective
algorithm of searching such points.

\textbf{3. }We say that the pyramid $ABCD$ is a Pythagorean-Diophantine
pyramid if the coordinates of the vertices are integers, the lengths of each
segment $AB,AC,BC,BD,AD,CD$ are natural numbers and the triangles $ADC,BDC$
and $ADB$ are Pythagorean. Are there Pythagorean-Diophantine pyramids?

\textbf{4. }We say that the quaternion $q=n+im+jr+ks,$ with $n,m,r,s$
$\in\mathbb{Z}$, is a Hamiltonian integer. Each $q$ is represented by a couple
of Gauss integers as follows: $q=z+wj$, where $z=n+im$ and $w=r+is.$ It is
interesting to examine the possibility to develop analogous theory for
Hamiltonian integers.

\textbf{5. }Examine the coloring problem for Diophantine carpets. For a large
class of such carpets the chromatic number is 2.

\subsubsection{ACKNOWLEDGMENTS}

The authors are grateful to R. Lazov for an useful discussion and the
suggestion of Problem 5 and to I.Tonov for his interest to this subject.

\subsection{REFERENCES}

[DT] S. Dimiev, I. Tonov, Diophantine Figures, Mathematics and Mathematical
Education, 15 (1986), 383-390 (in Bulgarian).

[B1] M.N.Brancheva, Diophantine Figures and Diophantine Carpets, Mathematics
and Mathematical Education, 30 (2001), 289-296.

[D] S.Dimiev, Physicalism in mathematics? Mathematics and Mathematical
Education, 30 (2001), 78-85.

[BDM] M. N. Brancheva, S. Dimiev, N. Milev, On Diophantine Figures, (submited
for publication)

[B2] M. N. Brancheva, Assymptotic for the function $\kappa(n)$ (submited for publication)

[IR] K. Ireland, M. Rosen, A Classical Introduction to ModernNumber Theory,
Springer Verlag, 1982

[TT] T. Tonkov, Figural Numbers (in bulgarian), Sofia, 1971

[P] K. Prachar, Primzahlverteilung, Springer Verlag, 1957

[E] H.M.Edwards, Fermat's Last Theorem, Springer Verlag, New-York, 1977

[G] P.Guncheva, Square Radicals of Gauss Integers (to be published).

\bigskip\lbrack DMY] S.Dimiev, K.Markov, M.Yawata, On the Geometry of
Pythagorean Triples (to be published)

\bigskip

Authors addresses:

S.Dimiev \ \ \ \ \ \ \ \ \ \ \ \ \ \ \ \ \ \ \ \ \ \ \ \ \ \ \ \ \ \ \ \ \ \ \ \ \ \ \ \ \ \ \ \ \ \ \ \ \ \ \ \ \ \ \ \ K.Markov

Institute of Mathematics and Informatics, \ \ \ \ \ \ \ \ \ Mladost-3, \ \ \ \ \ \ \ \ \ \ \ \ \ \ \ \ \ \ 

Bulgarian Academy of Sciences,
\ \ \ \ \ \ \ \ \ \ \ \ \ \ \ \ \ \ \ \ \ \ block 325, entry 6,

8, G.Bonchev str., Sofia 1113
\ \ \ \ \ \ \ \ \ \ \ \ \ \ \ \ \ \ \ \ \ \ \ \ \ \ Sofia 1117

e-mail: \ \ sdimiev@math.bas.bg
\ \ \ \ \ \ \ \ \ \ \ \ \ \ \ \ \ \ \ \ \ \ \ \ \ \ tel. (0357 2) 779133 \ \ \ 
\end{document}